\documentclass[twoside]{article}

\usepackage{epsfig}
\usepackage{amsfonts,amsmath,amsthm}
\usepackage[latin1]{inputenc}

\def\qed{\hbox to 0pt{}\hfill$\rlap{$\sqcap$}\sqcup$}

\newtheorem{thm}{Theorem}

\newcounter{subfigure}

\begin{document}

\vskip-45mm
\title{Gradual approximation of the domain of attraction by gradual extension of the
"embryo" of the transformed optimal Lyapunov function}

\author{E. Kaslik$^{1}$, A.M. Balint$^{2}$, \c{S}t. Balint$^{1}$ }

\date{}
\maketitle

\noindent$^1$Department of Mathematics,\hfill\break West
University of Timi\c{s}oara\hfill\break e-mail:
balint@balint.uvt.ro\hfill\break

\noindent$^2$Department of Physics,\hfill\break West University of
Timi\c{s}oara\hfill\break

\vspace{1mm}

\footnotetext{2000 \textit{AMS(MOS) Subject Classification:}
34.D.20; 34.D.45; 37.C.70}

\vskip 0.2in {\parindent=0pt \noindent{\bf Abstract} In this paper
an autonomous analytical system of ordinary differential equations
is considered. For an asymptotically stable steady state $x^{0}$
of the system a gradual approximation of the domain of attraction
($DA$) is presented in the case when the matrix of the linearized
system in $x^{0}$ is diagonalizable. This technique is based on
the gradual extension of the "embryo" of an analytic function of
several complex variables. The analytic function is the
transformed of a Lyapunov function whose natural domain of
analyticity is the $DA$ and which satisfies a linear
non-homogeneous partial differential equation. The equation
permits to establish an "embryo" of the transformed function and a
first approximation of $DA$. The "embryo" is used for the
determination of a new "embryo" and a new part of the $DA$. In
this way, computing new "embryos" and new domains, the $DA$ is
gradually approximated. Numerical examples are given for
polynomial systems. For systems considered recently in the
literature the results are compared with those obtained with other
methods.} \vskip 0.2in


\section{Introduction}
The domain of attraction ($DA$) of a steady state of a dynamical
system is the set of initial states from which the system
converges to the steady state itself. In order to guarantee stable
behavior of a dynamical system in a region of the state parameters
it is important to know the $DA$ [1].

Theoretical research shows that the $DA$ and its boundary are
complicated sets [2], [3], [4], [5]. In most cases, they do not
admit an explicit elementary representation. For this reason,
different procedures are used for the approximation of the $DA$
with domains having a simpler shape. This practice became
fundamental in the last 20 years [6]. The domain which
approximates the $DA$ is defined by a Lyapunov function, generally
quadratic. For a given Lyapunov function the computation of the
optimal estimate of the $DA$ amounts to solving a non-convex
distance problem [7], [8], [9], [10], [11], [12], [13].

In this paper we present a new technique for a gradual
approximation of the $DA$ in the case of autonomous analytical
systems of ordinary differential equations. More precisely for the
gradual approximation of the $DA$ of an asymptotically stable
steady state in which the matrix of the linearized system is
diagonalizable. This technique is based on the theoretical results
obtained in [14], [15], [16] and in the following we present a
summary of this theory.

\section{Theoretical summary}

We consider the system of differential equations

\begin{equation}
\label{ec}
\dot{x}=f(x)
\end{equation}
where $f:\mathbb{R}^{n}\rightarrow\mathbb{R}^{n}$ is an analytical
function with the following properties:
\begin{itemize}
    \item[i.] $f(0)=0$ (i.e. zero is a steady state of (\ref{ec}))
    \item[ii.] the real parts of the eigenvalues of $\frac{\partial f}{\partial x}(0)$
    are negative (i.e. $x=0$ is an asymptotically stable steady state)
\end{itemize}

The following theorem provides a tool of determining the $DA$ of
the zero steady state of (\ref{ec}).

\begin{thm}
(see [14]) The $DA$ of the null solution of (\ref{ec}) coincides
with the natural domain of analyticity of the unique solution V of
the problem

\begin{equation}
\label{partial}
\begin{array}{ll}
\left\{\begin{array}{l}
\langle\nabla V,f\rangle=-\|x\|^{2}\\
V(0)=0
\end{array}\right.
\end{array}
\end{equation}

The function $V$ is positive on $DA$ and
$\lim\limits_{x\rightarrow x_{0}}V(x)=\infty$ for any
$x_{0}\in\partial DA$.
\end{thm}

Thus, the problem of finding the $DA$ is reduced to the
determination of the natural domain of analyticity of the solution
$V$ of (\ref{partial}). This function will be called in the
followings the optimal Lyapunov function.

To determine the optimal Lyapunov function $V$ and its domain of
analyticity is not easy, but, in the diagonalizable case, we can
determine the coefficients of the expansion of $W=V\circ S$ in
$0$, where $S$ reduces $\frac{\partial f}{\partial x}(0)$ to the
diagonal form. Then, using a Cauchy-Hadamard type theorem, we can
obtain the domain of convergence $D^{0}$ of the series of $W$, and
$DA^{0}=S(D^{0})$ is a part of the domain of attraction.

\newpage

\subsection{The coefficients of the transformed optimal Lyapunov function in the diagonalizable case}

For the system (\ref{ec}) the following theorem holds:

\begin{thm}
(see [15]) For each isomorphism
$S:\mathbb{C}^{n}\rightarrow\mathbb{C}^{n}$ and $g=S^{-1}\circ
f\circ S$, the problem
\begin{equation}
\label{partial2}
\begin{array}{ll}
\left\{\begin{array}{l}
\langle\nabla W,g\rangle=-\|Sz\|^{2}\\
W(0)=0
\end{array}\right.
\end{array}
\end{equation}
has a unique analytical solution, namely $W=V\circ S$ where V is
the optimal Lyapunov function for (\ref{ec}).
\end{thm}

The function $W$ will be called the transformed optimal Lyapunov
function.

In the followings, we will suppose that the matrix $\frac{\partial
f}{\partial x}(0)$ is diagonalizable.

Let be $S:\mathbb{C}^{n}\rightarrow\mathbb{C}^{n}$ one isomorphism
which reduces $\frac{\partial f}{\partial x}(0)$ to the diagonal
form $S^{-1}\frac{\partial f}{\partial
x}(0)S=diag(\lambda_{1},\lambda_{2}...\lambda_{n})$ and
$g=S^{-1}\circ f\circ S$.

For this $S$, we consider the expansion of $W$ in $0$:

\begin{equation}\label{serieW}
    W(z_{1},z_{2},...,z_{n})=\sum\limits_{m=2}^{\infty}\sum\limits_{|j|=m}B_{j_{1}j_{2}...j_{n}}z_{1}^{j_{1}}z_{2}^{j_{2}}...z_{n}^{j_{n}}
\end{equation}

\noindent and the expansions in $0$ of the scalar components
$g_{i}$ of $g$:

\begin{equation}\label{serieg}
    g_{i}(z_{1},z_{2},...,z_{n})=\lambda_{i}y_{i}+\sum\limits_{m=2}^{\infty}\sum\limits_{|j|=m}b^{i}_{j_{1}j_{2}...j_{n}}z_{1}^{j_{1}}z_{2}^{j_{2}}...z_{n}^{j_{n}}
\end{equation}

\begin{thm}
(see [16]) The coefficients $B_{j_{1}j_{2}...j_{n}}$ of the
development (\ref{serieW}) are given by the following relations:
\begin{equation}
\label{coefB} B_{j_{1}j_{2}...j_{n}}=
\begin{array}{lll}
\left\{\begin{array}{l}
-\frac{1}{2\lambda_{i_{0}}}\sum\limits_{i=1}^{n}s^{2}_{ii_{0}}
\textrm{
if   } |j|=j_{i_{0}}=2\\ \\
-\frac{2}{\lambda_{p}+\lambda_{q}}\sum\limits_{i=1}^{n}s_{ip}s_{iq}
\textrm{  if   } |j|=2 \textrm{ and } j_{p}=j_{q}=1\\ \\
-\frac{1}{\sum\limits_{i=1}^{n}j_{i}\lambda_{i}}\sum\limits_{p=2}^{|j|-1}\sum\limits_{|k|=p,k_{i}\leq
j_{i}}\sum\limits_{i=1}^{n}[(j_{i}-k_{i}+1)\\
b^{i}_{j_{1}j_{2}...j_{n}}B_{j_{1}-k_{1}...j_{i}-k_{i}+1...j_{n}-k_{n}}]
\textrm{ if } |j|\geq3
\end{array}\right.
\end{array}
\end{equation}
\end{thm}

According to (\ref{coefB}), the coefficients of the terms of
degree $m\geq 3$ are obtained in function of the coefficients of
the terms of degree smaller than $m$.

\newpage

\subsection{The domain of convergence of the series of W}

We define the domain of convergence of the series (\ref{serieW})
as the set of those $z\in\mathbb{C}^{n}$ which have the property
that the series (\ref{serieW}) is absolutely convergent in a
neighborhood of $z$ [17]. We denote by $D^{0}$ this domain.
Actually, the domain of convergence $D^{0}$ coincides with the
interior of the set of all the points $z^{0}$ in which the series
(\ref{serieW}) is absolutely convergent (see [17]).

\begin{thm}
(Cauchy-Hadamard see [17]) A point $z$ belongs to $D^{0}$ if and
only if
\begin{equation} \label{ineg}
\overline{\lim_{m}}\sqrt[m]{\sum_{|j|=m}|B_{j_{1}j_{2}...j_{n}}z_{1}^{j_{1}}z_{2}^{j_{2}}...z_{n}^{j_{n}}|}<1
\end{equation}
\end{thm}

Using $D^{0}$ we can find a part of the $DA$.

\begin{thm}
(see [15]) If $z$ belongs to $D^{0}$ then the point $x=Sz$ is in
the domain of attraction $DA$.
\end{thm}

Actually, by the previous theorem, we obtain the subdomain
$S(D^{0})\subset DA$, which has the property that $\partial
S(D^{0})\cap\partial DA\neq\emptyset$. We also observe that the
subdomain $S(D^{0})$ is symmetrical to the origin (actually, the
domain of convergence $D^{0}$ is symmetrical to the origin and the
axes of coordinates). This first approximation of the domain of
attraction $DA$ will be denoted by $DA^{0}=S(D^{0})$.

Condition (\ref{ineg}) represents an algoritmizable criterion for
determining the region of convergence $D^{0}$.

\subsection{A method of extending the estimate of DA}

In practice, we can compute the coefficients
$B_{j_{1}j_{2}...j_{n}}$ up to a finite degree $p$. This degree
$p$ has to be big enough for assuming that the domain $D^{0}_{p}$
given by

\begin{equation}
D^{0}_{p}=\{z\in\mathbb{C}^{n}/\sqrt[p]{\sum\limits_{|j|=p}|B_{j_{1}j_{2}...j_{n}}
 z_{1}^{j_{1}}z_{2}^{j_{2}}...z_{n}^{j_{n}}|}<1\}
 \nonumber
\end{equation}

\noindent approximates the region of convergence $D^{0}$ of the
series of $W$. On this domain $D^{0}_{p}$ we approximate $W$ by
the "embryo"

\begin{equation}
 W^{0}_{p}(z_{1},z_{2},...,z_{n})=\sum\limits_{m=2}^{p}\sum\limits_{|j|=m}B_{j_{1}j_{2}...j_{n}}
 z_{1}^{j_{1}}z_{2}^{j_{2}}...z_{n}^{j_{n}}
\end{equation}

The first estimate of the region of attraction $DA$ will be
$DA^{0}_{p}=S(D^{0}_{p})$.

In order to extend the first estimate $DA^{0}_{p}$, we will expand
$W$ in a point $z^{0}$ close to $\partial D^{0}_{p}$ in which
$|W^{0}_{p}(z^{0})|$ is still small. That is because, according to
Theorem 1, the points $z$ close to $\partial D^{0}_{p}$ for which
$|W^{0}_{p}(z)|$ is extremely high are close to $\partial
S^{-1}(DA)$.

To find the expansion of $W$ in $z^{0}$ close to $\partial
D^{0}_{p}$, we will compute the expansion in $z^{0}$ of the
"embryo" $W^{0}_{p}$ of $W$. We obtain:

\begin{eqnarray}
W^{1}_{p}(z_{1},z_{2},...,z_{n})\!\!\!\!\!\!&=&\!\!\!\!\!\!
\sum\limits_{m=0}^{p}\sum\limits_{|j|=m}\frac{\partial^{m}W_{p}}{\partial
z_{1}^{j_{1}}\partial z_{2}^{j_{2}}...\partial
z_{n}^{j_{n}}}(z^{0})
 (z_{1}-z_{1}^{0})^{j_{1}}(z_{2}-z_{2}^{0})^{j_{2}}...(z_{n}-z_{n}^{0})^{j_{n}}
 /m!=\nonumber\\
 \!\!\!\!\!\!&=&\!\!\!\!\!\!
 \sum\limits_{m=0}^{p}\sum\limits_{|j|=m}B^{1}_{j_{1}j_{2}...j_{n}}
 (z_{1}-z_{1}^{0})^{j_{1}}(z_{2}-z_{2}^{0})^{j_{2}}...(z_{n}-z_{n}^{0})^{j_{n}}
\end{eqnarray}

\noindent We consider the set

\begin{equation}
D_{p}^{1}=\{z\in\mathbb{C}^{n}/\sqrt[p]{\sum\limits_{|j|=p}|B_{j_{1}j_{2}...j_{n}}
 (z_{1}-z_{1}^{0})^{j_{1}}(z_{2}-z_{2}^{0})^{j_{2}}...(z_{n}-z_{n}^{0})^{j_{n}}|}<1\}
 \nonumber
\end{equation}

\noindent which provides a new part $DA^{1}_{p}=S(D_{p}^{1})$ of
$DA$.

So, $DA_{p}^{0}\cup DA_{p}^{1}$ gives a larger estimate of the
$DA$. We can continue this procedure for a few steps, till the
values $|W_{p}^{k}|$ become extremely large and we obtain the
estimate $DA_{p}^{0}\cup DA_{p}^{1}\cup ... \cup
DA_{p}^{k}=S(D^{0}_{p}\cup D_{p}^{1}\cup ... \cup D_{p}^{k})$ of
$DA$.

\section{Numerical results}

The computations were made using a program written in Mathematica
4, Wolfram Research. The data for the estimations (the degree up
to which the approximation is made, the necessary timing for the
estimations) are displayed in Table~1.

\subsection{Systems with known domains of attraction}

In this subsection, we will present some examples of systems of
two or three differential equations, for which we can compute
easily the $DA$. We will apply our technique to these examples,
and we will show how the real domains of attraction are gradually
approximated. These examples are meant to validate our procedure.

\subsubsection{Example 1}

This is an example of a system for which the null solution has a
bounded domain of attraction:

\begin{equation}
\label{ex2}
\begin{array}{ll}
\left\{\begin{array}{l}
\dot{x_{1}}=-x_{1}[4-(x_{1}-1)^{2}-x_{2}^{2}]\\
\dot{x_{2}}=-x_{2}[4-(x_{1}-1)^{2}-x_{2}^{2}]
\end{array}\right.
\end{array}
\end{equation}

\noindent The domain of attraction of the null solution of this
system is the interior of the circle of radius $2$ centered in
$(1,0)$:
\begin{displaymath}
DA=\{(x_{1},x_{2})\in\mathbb{R}^{2}/(x_{1}-1)^{2}+x_{2}^{2}<4\}
\end{displaymath}

\noindent After three steps, we obtain the estimate shown in
Figure 1.

\newpage

The thick black line represents the true boundary of the domain of
attraction, the dark grey set denotes the first estimate
$DA_{p}^{0}$ of $DA$ and the further estimates $DA_{p}^{k}$ of
$DA$ with $k\geq 1$ are colored in light grey.

\renewcommand{\thefigure}{\arabic{figure}}

\setcounter{subfigure}{1}

\vspace*{-0.5cm}

\begin{figure}[htbp]
\centering \includegraphics*[bb=1.5cm 0cm 8.cm
9cm,width=4.8cm,height=4.8cm,angle=0]{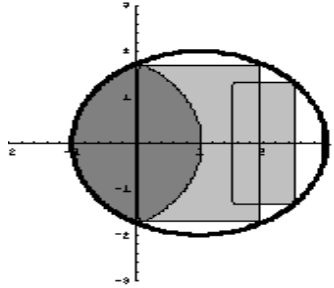} \caption{The
estimate of $DA$ obtained after three steps for
system~(\ref{ex2})}
\end{figure}

\subsubsection{Example 2}

The following system of three differential equations is
considered:

\begin{equation}
\label{ex4}
\begin{array}{ll}
\left\{\begin{array}{l}
\dot{x_{1}}=-x_{1}(1-x_{1}^{2}-x_{2}^{2}+x_{3}^{2})\\
\dot{x_{2}}=-x_{2}(1-x_{1}^{2}-x_{2}^{2}+x_{3}^{2})\\
\dot{x_{3}}=-x_{3}(1-x_{1}^{2}-x_{2}^{2}+x_{3}^{2})
\end{array}\right.
\end{array}
\end{equation}

\noindent The boundary of the $DA$ of the null solution of this
system is:
\begin{displaymath}
\partial
DA=\{(x_{1},x_{2},x_{3})\in\mathbb{R}^{3}/1-x_{1}^{2}-x_{2}^{2}+x_{3}^{2}=0\}
\end{displaymath}

\noindent The first estimate $DA_{p}^{0}$ is shown in Figure 2.1.
After 2 steps, we obtain the estimate shown in Figure 2.2.

\renewcommand{\thefigure}{\arabic{figure}.\arabic{subfigure}}
\setcounter{subfigure}{1}

\begin{figure}[htbp]
\begin{minipage}[t]{0.5\linewidth}
\centering \includegraphics*[bb=0cm 0cm 9.5cm
10.5cm,width=4.8cm,height=4.8cm,angle=0]{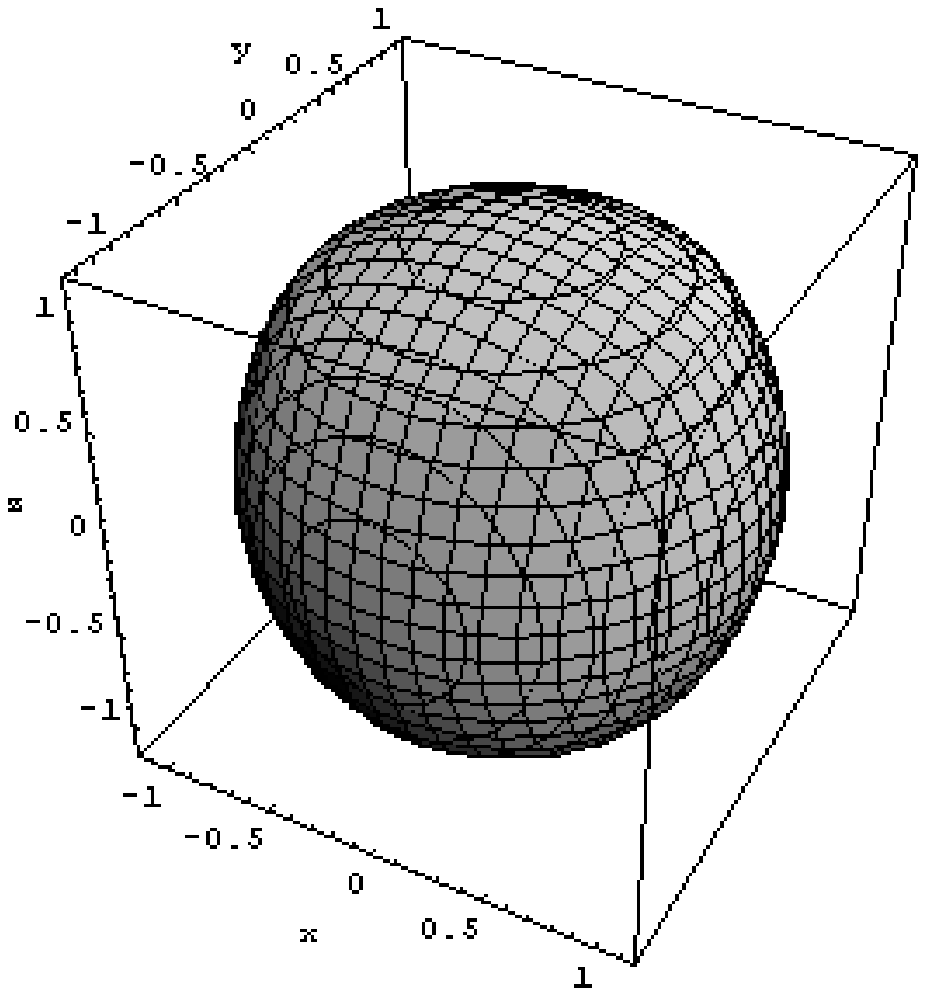}\caption{The
estimate of $DA$ obtained after 1 step for system~(\ref{ex4})}
\end{minipage}
\addtocounter{figure}{-1} \addtocounter{subfigure}{1}
\begin{minipage}[t]{0.5\linewidth}
\centering \includegraphics*[bb=0cm 0cm 8cm
10.5cm,width=4.8cm,height=4.8cm,angle=0]{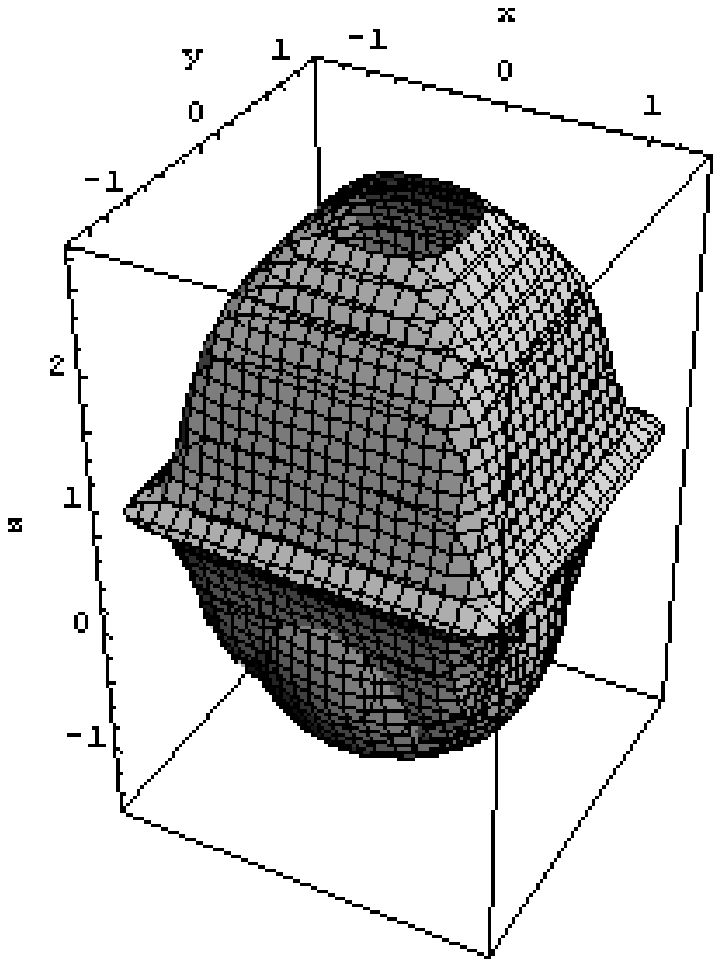}\caption{The
estimate of $DA$ obtained after 2 steps for system~(\ref{ex4})}
\end{minipage}
\end{figure}

\vspace*{-1cm}

\newpage

\subsection{Systems for which we don't know the DA}

In this subsection, some systems of differential equations are
presented for which we don't know the $DA$. For these examples, we
will apply the technique presented in Section 2, and we will show
that the estimate obtained using this technique is better than the
estimates obtained in [10].

\subsubsection{Example 3}

In [10], the following example is considered:

\begin{equation}
\label{ex5}
\begin{array}{ll}
\left\{\begin{array}{l} \dot{x_{1}}=x_{2}\\
\dot{x_{2}}=-2x_{1}-3x_{2}+x_{1}^{2}-x_{2}^2+x_{1}x_{2}
\end{array}\right.
\end{array}
\end{equation}

\noindent In Figure 3.1, an estimate of the $DA$ is shown obtained
after two steps for three different points close to the boundary
of $DA_{p}^{0}$. We observe that this estimate covers the estimate
presented in [10]. Figure 3.2 presents the estimate of $DA$
obtained after four steps. The thick black lines plotted in the
following figures represent a part of the approximated boundary of
the $DA$. The black interrupted line represents the boundary of
the estimate of the domain of attraction obtained in [10].

\setcounter{subfigure}{1}

\begin{figure}[htbp]
\begin{minipage}[t]{0.5\linewidth}
\centering \includegraphics*[bb=0cm 0cm 9cm
10cm,width=4.8cm,height=4.8cm,angle=0]{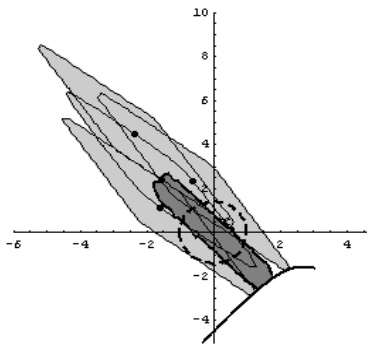} \caption{The
estimate of $DA$ after 2 steps for three different points close to
the boundary of $DA_{p}^{0}$ for system~(\ref{ex5})}
\end{minipage}
\addtocounter{figure}{-1} \addtocounter{subfigure}{1}
\begin{minipage}[t]{0.5\linewidth}
\centering \includegraphics*[bb=0cm 0cm 8.5cm
9cm,width=4.8cm,height=4.8cm,angle=0]{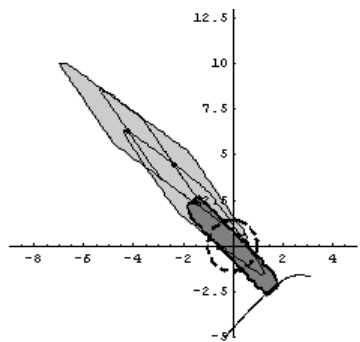}\caption{The
estimate of $DA$ obtained after 4 steps for system~(\ref{ex5})}
\end{minipage}
\end{figure}

\subsubsection{Example 4}

In [10], the following system of three equations is considered:

\begin{equation}
\label{ex8}
\begin{array}{ll}
\left\{\begin{array}{l}
\dot{x_{1}}=x_{2}\\
\dot{x_{2}}=x_{3}\\
\dot{x_{3}}=-3x_{1}-3x_{2}-2x_{3}+x_{1}^{3}+x_{2}^{3}+x_{3}^{3}
\end{array}\right.
\end{array}
\end{equation}

\noindent Figure 4.1 shows the estimate of $DA$ obtained after one
step. The estimate of the $DA$ obtained after two steps is
presented in Fig. 4.2.

\setcounter{subfigure}{1}

\begin{figure}[htbp]
\begin{minipage}[t]{0.5\linewidth}
\centering \includegraphics*[bb=0cm 0cm 7cm
11cm,width=4.8cm,height=4.8cm,angle=0]{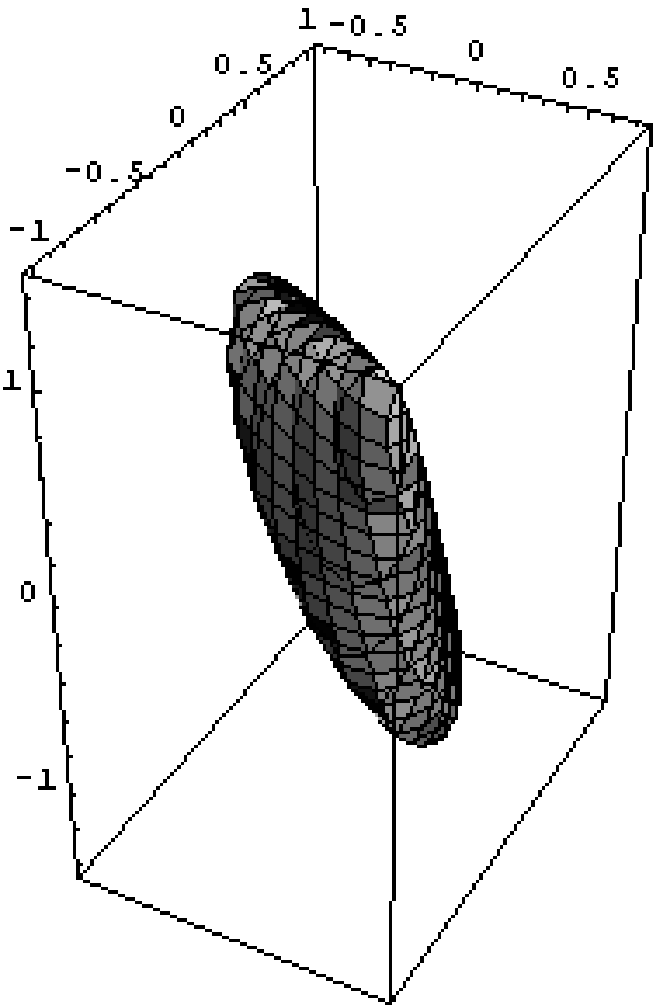}\caption{The
estimate of $DA$ obtained after 1 step for system~(\ref{ex8})}
\end{minipage}
\addtocounter{figure}{-1} \addtocounter{subfigure}{1}
\begin{minipage}[t]{0.5\linewidth}
\centering \includegraphics*[bb=0cm 0cm 7cm
12cm,width=4.8cm,height=4.8cm,angle=0]{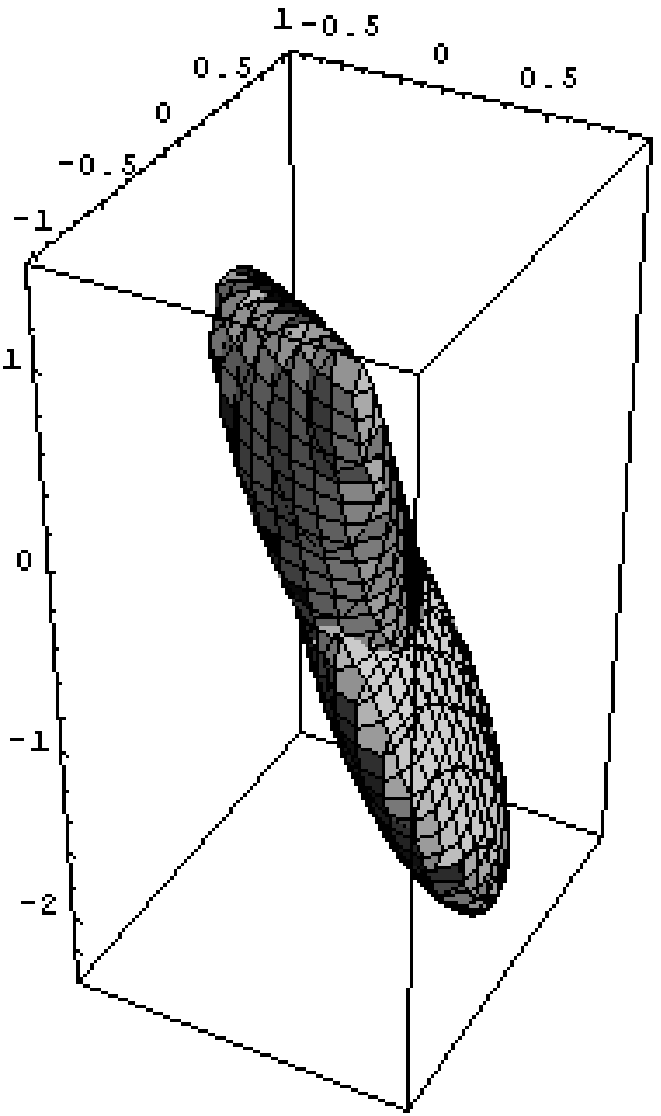}\caption{The
estimate of $DA$ obtained after 2 steps for system~(\ref{ex8})}
\end{minipage}
\end{figure}

\newpage

\begin{center}
Table 1. Numerical data
\end{center}

\begin{tabular}{|c|c|c|c|}
  \hline
  example & order of approximation & timing for $1^{st}$ step & timing for $2^{nd}$ step \\
  \hline
  \hline
  1 & 50 & 1.7 min & 34.7 min \\
  \hline
  2 & 30 & 10.2 min & 14.7 h \\
  \hline
  3 & 30 & 0.9 min & 9.9 min \\
 \hline
  4 & 30 & 19.1 min & 16.2 h \\
  \hline
\end{tabular}

\end{document}